\documentclass[10pt,draft]{article}
\usepackage{amsfonts,amsmath,amssymb,cite}
\numberwithin{equation}{section}
\begin{document}
\title{Some differentiation formulas for Legendre polynomials}
\author{Rados{\l}aw Szmytkowski \\*[3ex]
Atomic Physics Division, 
Department of Atomic Physics and Luminescence, \\
Faculty of Applied Physics and Mathematics,
Gda{\'n}sk University of Technology, \\
Narutowicza 11/12, PL 80--233 Gda{\'n}sk, Poland \\
email: radek@mif.pg.gda.pl}
\date{\today}
\maketitle
\begin{abstract}
In a series of recent works, we have provided a number of explicit
expressions for the derivative of the associated Legendre function of
the first kind with respect to its degree, 
$[\partial P_{\nu}^{m}(z)/\partial\nu]_{\nu=n}$, with
$m,n\in\mathbb{N}$. In this communication, we use some of those
expressions to obtain several, we believe new, explicit formulas for
the derivatives $\mathrm{d}^{m}[P_{n}(z)\ln(z\pm1)]/\mathrm{d}z^{m}$,
where $P_{n}(z)$ is the Legendre polynomial.
\vskip1ex
\noindent
\textbf{KEY WORDS:} Legendre polynomials; Legendre functions; 
special functions
\vskip1ex
\noindent
\textbf{MSC2010:} 33C45, 33C55
\end{abstract}
\maketitle
%
%
\section{The problem and the method}
\label{I}
\setcounter{equation}{0}
Recently, Brychkov has published a monumental reference work
\cite{Bryc08} containing, among others, a large number of explicit
expressions for derivatives of various special functions with respect
to their arguments and parameters. It is the purpose of the present
communication to supplement the handbook \cite{Bryc08} with several
closed-form formulas for the derivatives
$\mathrm{d}^{m}[P_{n}(z)\ln(z\pm1)]/\mathrm{d}z^{m}$, where
$P_{n}(z)$ is the Legendre polynomial of degree $n$. We shall arrive
at these formulas exploiting results of recent papers
\cite{Szmy06,Szmy07,Szmy09a,Szmy09b}, in which we have extensively
investigated the derivatives $[\partial
P_{\nu}(z)/\partial\nu]_{\nu=n}$ and $[\partial
P_{\nu}^{m}(z)/\partial\nu]_{\nu=n}$, where $P_{\nu}(z)$ and
$P_{\nu}^{m}(z)$ are the Legendre function of the first kind and the
associated Legendre function of the first kind, respectively.

To make the functions which appear in the following considerations
single-valued, we cut the complex $z$-plane along the real axis from
$-\infty$ to $+1$. Then, it follows in particular that
\begin{equation}
-z-1=\mathrm{e}^{\mp\mathrm{i}\pi}(z+1),
\qquad
-z+1=\mathrm{e}^{\mp\mathrm{i}\pi}(z-1)
\qquad (\arg(z)\gtrless0).
\label{1.1}
\end{equation}
Throughout the paper, it is assumed that $\nu\in\mathbb{C}$ and
$m,n\in\mathbb{N}$.

We begin with recalling the following formulas
\cite{Szmy06,Szmy07,Szmy09a,Szmy09b}:
\begin{equation}
\frac{\partial P_{\nu}(z)}{\partial\nu}\bigg|_{\nu=n}
=P_{n}(z)\ln\frac{z+1}{2}+R_{n}(z)
\label{1.2}
\end{equation}
and
\begin{equation}
\frac{\partial P_{\nu}^{m}(z)}{\partial\nu}\bigg|_{\nu=n}
=P_{n}^{m}(z)\ln\frac{z+1}{2}+R_{n}^{m}(z)
\label{1.3}
\end{equation}
for the derivatives of $P_{\nu}(z)$ and $P_{\nu}^{m}(z)$ with respect
to their degrees $\nu$, the two functions being related through
\begin{equation}
P_{\nu}^{m}(z)=(z^{2}-1)^{m/2}
\frac{\mathrm{d}^{m}P_{\nu}(z)}{\mathrm{d}z^{m}}.
\label{1.4}
\end{equation}
The relationship analogous to that in Eq.\ (\ref{1.4}) holds also
between the function $P_{n}^{m}(z)$ and the polynomial $P_{n}(z)$. In
Eqs.\ (\ref{1.2}) and (\ref{1.3}), $R_{n}(z)$ is a known polynomial
in $z$ of degree $n$ (the Bromwich polynomial) and $R_{n}^{m}(z)$ is
a known function (see Section \ref{III} below). It holds that
\begin{equation}
R_{n}^{0}(z)=R_{n}(z),
\label{1.5}
\end{equation}
but it must be emphasized that in general $R_{n}^{m}(z)$ is
\emph{not\/} related to $R_{n}(z)$ through a formula analogous to
that in Eq.\ (\ref{1.4}).

Differentiating Eq.\ (\ref{1.2}) $m$ times with respect to $z$ gives
\begin{equation}
\frac{\mathrm{d}^{m}}{\mathrm{d}z^{m}}
\frac{\partial P_{\nu}(z)}{\partial\nu}\bigg|_{\nu=n}
=\frac{\mathrm{d}^{m}}{\mathrm{d}z^{m}}
\left[P_{n}(z)\ln\frac{z+1}{2}\right]
+\frac{\mathrm{d}^{m}R_{n}(z)}{\mathrm{d}z^{m}}.
\label{1.6}
\end{equation}
On the other hand, from Eqs.\ (\ref{1.3}) and (\ref{1.4}) it follows
that
\begin{equation}
\frac{\mathrm{d}^{m}}{\mathrm{d}z^{m}}
\frac{\partial P_{\nu}(z)}{\partial\nu}\bigg|_{\nu=n}
=\frac{\mathrm{d}^{m}P_{n}(z)}{\mathrm{d}z^{m}}\ln\frac{z+1}{2}
+(z^{2}-1)^{-m/2}R_{n}^{m}(z).
\label{1.7}
\end{equation}
Combining Eqs.\ (\ref{1.6}) and (\ref{1.7}) gives
\begin{equation}
\frac{\mathrm{d}^{m}}{\mathrm{d}z^{m}}[P_{n}(z)\ln(z+1)]
=\frac{\mathrm{d}^{m}P_{n}(z)}{\mathrm{d}z^{m}}\ln(z+1)
+(z^{2}-1)^{-m/2}R_{n}^{m}(z)
-\frac{\mathrm{d}^{m}R_{n}(z)}{\mathrm{d}z^{m}}.
\label{1.8}
\end{equation}
If $m>n$, the two derivatives on the right-hand side of Eq.\
(\ref{1.8}) vanish and one simply has
\begin{equation}
\frac{\mathrm{d}^{m}}{\mathrm{d}z^{m}}[P_{n}(z)\ln(z+1)]
=(z^{2}-1)^{-m/2}R_{n}^{m}(z)
\qquad (m>n).
\label{1.9}
\end{equation}
Replacement of $z$ by $-z$ in Eqs.\ (\ref{1.8}) and (\ref{1.9}),
followed by the use of Eq.\ (\ref{1.1}) and the well-known property
\begin{equation}
P_{n}(-z)=(-)^{n}P_{n}(z),
\label{1.10}
\end{equation}
implies the counterpart relationships
\begin{equation}
\frac{\mathrm{d}^{m}}{\mathrm{d}z^{m}}[P_{n}(z)\ln(z-1)]
=\frac{\mathrm{d}^{m}P_{n}(z)}{\mathrm{d}z^{m}}\ln(z-1)
+(-)^{n}(z^{2}-1)^{-m/2}R_{n}^{m}(-z)
-(-)^{n}\frac{\mathrm{d}^{m}R_{n}(-z)}{\mathrm{d}z^{m}}
\label{1.11}
\end{equation}
and
\begin{equation}
\frac{\mathrm{d}^{m}}{\mathrm{d}z^{m}}[P_{n}(z)\ln(z-1)]
=(-)^{n}(z^{2}-1)^{-m/2}R_{n}^{m}(-z)
\qquad (m>n).
\label{1.12}
\end{equation}
Thus, from Eqs.\ (\ref{1.8}), (\ref{1.9}), (\ref{1.11}) and
(\ref{1.12}) we see that once the function $R_{n}^{m}(z)$ and the
polynomial $R_{n}(z)$ are known, the derivatives
$\mathrm{d}^{m}[P_{n}(z)\ln(z\pm1)]/\mathrm{d}z^{m}$ may be
evaluated.
\section{Explicit representations of the polynomial $R_{n}(z)$ and
the function $R_{n}^{m}(z)$} 
\label{II} 
\setcounter{equation}{0}
The following three representations of the Bromwich polynomial
$R_{n}(z)$ have been derived in Refs.\ \cite{Szmy06,Szmy07} (in the
first of these papers, the reader will find coordinates of earlier
publications of Schelkunoff and Bromwich, in which the expressions
(\ref{2.1}) and (\ref{2.3}) were found differently than in Ref.\
\cite{Szmy06}):
\begin{equation}
R_{n}(z)=-2\psi(n+1)P_{n}(z)
+2\sum_{k=0}^{n}\frac{(k+n)!\psi(k+n+1)}{(k!)^{2}(n-k)!}
\left(\frac{z-1}{2}\right)^{k},
\label{2.1}
\end{equation}
\begin{equation}
R_{n}(z)=2\sum_{k=0}^{n}(-)^{k+n}\frac{(k+n)!}{(k!)^{2}(n-k)!}
[\psi(k+n+1)-\psi(k+1)]\left(\frac{z+1}{2}\right)^{k},
\label{2.2}
\end{equation}
\begin{equation}
R_{n}(z)=2[\psi(2n+1)-\psi(n+1)]P_{n}(z)
+2\sum_{k=0}^{n-1}(-)^{k+n}\frac{2k+1}{(n-k)(k+n+1)}P_{k}(z),
\label{2.3}
\end{equation}
where
\begin{equation}
\psi(\zeta)=\frac{1}{\Gamma(\zeta)}
\frac{\mathrm{d}\Gamma(\zeta)}{\mathrm{d}\zeta}
\label{2.4}
\end{equation}
is the digamma function. Next, in Ref.\ \cite{Szmy09a} it has been
proved that the function $R_{n}^{m}(z)$ may be written as
\begin{eqnarray}
R_{n}^{m}(z) &=& -\,[\psi(n+1)+\psi(n-m+1)]P_{n}^{m}(z)
\nonumber \\
&& +\,\left(\frac{z^{2}-1}{4}\right)^{m/2}
\sum_{k=0}^{n-m}\frac{(k+n+m)!\psi(k+n+m+1)}{k!(k+m)!(n-m-k)!}
\left(\frac{z-1}{2}\right)^{k}
\nonumber \\
&& +\,\frac{(n+m)!}{(n-m)!}\left(\frac{z-1}{z+1}\right)^{m/2}
\sum_{k=0}^{n}\frac{(k+n)!\psi(k+n+1)}{k!(k+m)!(n-k)!}
\left(\frac{z-1}{2}\right)^{k}
\nonumber \\
&& (0\leqslant m\leqslant n),
\label{2.5}
\end{eqnarray}
while in Ref.\ \cite{Szmy09b} the following two alternative
expressions:
\begin{eqnarray}
R_{n}^{m}(z) &=& -\,[\psi(n+m+1)+\psi(n-m+1)]P_{n}^{m}(z)
\nonumber \\
&& +\,\left(\frac{z^{2}-1}{4}\right)^{m/2}
\sum_{k=0}^{n-m}\frac{(k+n+m)!}{k!(k+m)!(n-m-k)!}
\nonumber \\
&& \quad
\times[2\psi(k+n+m+1)-\psi(k+m+1)]\left(\frac{z-1}{2}\right)^{k}
\nonumber \\
&& +\,\frac{(n+m)!}{(n-m)!}\left(\frac{z-1}{z+1}\right)^{m/2}
\sum_{k=0}^{n}\frac{(k+n)!\psi(k+m+1)}{k!(k+m)!(n-k)!}
\left(\frac{z-1}{2}\right)^{k}
\nonumber \\
&& (0\leqslant m\leqslant n)
\label{2.6}
\end{eqnarray}
and
\begin{eqnarray}
R_{n}^{m}(z) &=& [\psi(n+m+1)-2\psi(n+1)-\psi(n-m+1)]P_{n}^{m}(z)
\nonumber \\
&& +\,\left(\frac{z^{2}-1}{4}\right)^{m/2}
\sum_{k=0}^{n-m}\frac{(k+n+m)!\psi(k+m+1)}{k!(k+m)!(n-m-k)!}
\left(\frac{z-1}{2}\right)^{k}
\nonumber \\
&& +\,\frac{(n+m)!}{(n-m)!}\left(\frac{z-1}{z+1}\right)^{m/2}
\sum_{k=0}^{n}\frac{(k+n)!}{k!(k+m)!(n-k)!}
\nonumber \\
&& \quad
\times[2\psi(k+n+1)-\psi(k+m+1)]\left(\frac{z-1}{2}\right)^{k}
\qquad (0\leqslant m\leqslant n)
\label{2.7}
\end{eqnarray}
have been provided. Furthermore, in Ref.\ \cite{Szmy09b} we have
found the formula
\begin{eqnarray}
R_{n}^{m}(z) &=& 
-\,(-)^{n+m}\left(\frac{z-1}{z+1}\right)^{m/2}\sum_{k=0}^{m-1}
\frac{(k+n)!(m-k-1)!}{k!(n-k)!}\left(\frac{z+1}{2}\right)^{k}
\nonumber \\
&& +\,(-)^{n+m}\left(\frac{z^{2}-1}{4}\right)^{m/2}
\sum_{k=0}^{n-m}(-)^{k}\frac{(k+n+m)!}{k!(k+m)!(n-m-k)!}
\nonumber \\
&& \quad \times[2\psi(k+n+m+1)-\psi(k+m+1)-\psi(k+1)]
\left(\frac{z+1}{2}\right)^{k}
\nonumber \\
&& (0\leqslant m\leqslant n),
\label{2.8}
\end{eqnarray}
while in Ref.\ \cite{Szmy09a} we have arrived at
\begin{eqnarray}
R_{n}^{m}(z) &=& [2\psi(2n+1)-\psi(n+1)-\psi(n-m+1)]P_{n}^{m}(z)
\nonumber \\
&& +\,(-)^{n}\frac{(n+m)!}{(n-m)!}\sum_{k=0}^{m-1}(-)^{k}
\frac{2k+1}{(n-k)(k+n+1)}P_{k}^{-m}(z)
\nonumber \\
&& +\,(-)^{n+m}\sum_{k=0}^{n-m-1}(-)^{k}
\frac{2k+2m+1}{(n-m-k)(k+n+m+1)}
\nonumber \\
&& \quad \times\left[1+\frac{k!(n+m)!}{(k+2m)!(n-m)!}\right]
P_{k+m}^{m}(z)
\qquad (0\leqslant m\leqslant n).
\label{2.9}
\end{eqnarray}
Finally, in Ref.\ \cite{Szmy09a} we have also obtained the remarkably
simple representation
\begin{equation}
R_{n}^{m}(z)=(-)^{n+m+1}(n+m)!(m-n-1)!P_{n}^{-m}(z)
\qquad (m>n).
\label{2.10}
\end{equation}
It is easy to see that for $m=0$ the triple of equations
(\ref{2.5})--(\ref{2.7}) reduces to Eq.\ (\ref{2.1}), while Eqs.\
(\ref{2.8}) and (\ref{2.9}) go over into Eqs.\ (\ref{2.2}) and
(\ref{2.3}), respectively, in accordance with Eq.\ (\ref{1.5}).
\section{Closed-form expressions for the derivatives \\
$\mathrm{d}^{m}[P_{n}(z)\ln(z\pm1)]/\mathrm{d}z^{m}$}
\label{III}
\setcounter{equation}{0}
Inserting the expansions (\ref{2.1})--(\ref{2.3}) and
(\ref{2.5})--(\ref{2.10}) into Eqs.\ (\ref{1.8}), (\ref{1.9}),
(\ref{1.11}) and (\ref{1.12}), and using the relation \cite[Eq.\
(8.936.2)]{Grad94}
\begin{equation}
\frac{\mathrm{d}^{m}P_{n}(z)}{\mathrm{d}z^{m}}
=\frac{(2m)!}{2^{m}m!}C_{n-m}^{(m+1/2)}(z),
\label{3.1}
\end{equation}
where $C_{n}^{(\alpha)}(z)$ is the Gegenbauer polynomial, yields the
following representations of the derivatives $\mathrm{d}^{m}
[P_{n}(z)\ln(z\pm1)]/\mathrm{d}z^{m}$:
\begin{eqnarray}
\frac{\mathrm{d}^{m}}{\mathrm{d}z^{m}}[P_{n}(z)\ln(z\pm1)]
&=& \frac{(2m)!}{2^{m}m!}C_{n-m}^{(m+1/2)}(z)\ln(z\pm1)
\nonumber \\
&& +\,\frac{(2m)!}{2^{m}m!}[\psi(n+1)-\psi(n-m+1)]
C_{n-m}^{(m+1/2)}(z)
\nonumber \\
&& -\,\frac{(\pm)^{n+m}}{2^{m}}
\sum_{k=0}^{n-m}(\pm)^{k}\frac{(k+n+m)!\psi(k+n+m+1)}
{k!(k+m)!(n-m-k)!}\left(\frac{z\mp1}{2}\right)^{k}
\nonumber \\
&& +\,(\pm)^{n}\,\frac{(n+m)!}{(n-m)!}(z\pm1)^{-m}
\sum_{k=0}^{n}(\pm)^{k}\frac{(k+n)!\psi(k+n+1)}
{k!(k+m)!(n-k)!}\left(\frac{z\mp1}{2}\right)^{k}
\nonumber \\
&& (0\leqslant m\leqslant n),
\label{3.2}
\end{eqnarray}
\begin{eqnarray}
\frac{\mathrm{d}^{m}}{\mathrm{d}z^{m}}[P_{n}(z)\ln(z\pm1)]
&=& \frac{(2m)!}{2^{m}m!}C_{n-m}^{(m+1/2)}(z)\ln(z\pm1)
\nonumber \\
&& -\,\frac{(2m)!}{2^{m}m!}[\psi(n+m+1)-2\psi(n+1)+\psi(n-m+1)]
C_{n-m}^{(m+1/2)}(z)
\nonumber \\
&& -\,\frac{(\pm)^{n+m}}{2^{m}}
\sum_{k=0}^{n-m}(\pm)^{k}\frac{(k+n+m)!\psi(k+m+1)}
{k!(k+m)!(n-m-k)!}\left(\frac{z\mp1}{2}\right)^{k}
\nonumber \\
&& +\,(\pm)^{n}\frac{(n+m)!}{(n-m)!}(z\pm1)^{-m}
\sum_{k=0}^{n}(\pm)^{k}\frac{(k+n)!\psi(k+m+1)}
{k!(k+m)!(n-k)!}\left(\frac{z\mp1}{2}\right)^{k}
\nonumber \\
&& (0\leqslant m\leqslant n),
\label{3.3}
\end{eqnarray}
\begin{eqnarray}
\frac{\mathrm{d}^{m}}{\mathrm{d}z^{m}}[P_{n}(z)\ln(z\pm1)]
&=& \frac{(2m)!}{2^{m}m!}C_{n-m}^{(m+1/2)}(z)\ln(z\pm1)
\nonumber \\
&& +\,\frac{(2m)!}{2^{m}m!}[\psi(n+m+1)-\psi(n-m+1)]
C_{n-m}^{(m+1/2)}(z)
\nonumber \\
&& -\,\frac{(\pm)^{n+m}}{2^{m}}\sum_{k=0}^{n-m}(\pm)^{k}
\frac{(k+n+m)!}{k!(k+m)!(n-m-k)!}
\nonumber \\
&& \quad \times[2\psi(k+n+m+1)-\psi(k+m+1)]
\left(\frac{z\mp1}{2}\right)^{k}
\nonumber \\
&& +\,(\pm)^{n}\frac{(n+m)!}{(n-m)!}(z\pm1)^{-m}
\sum_{k=0}^{n}(\pm)^{k}\frac{(k+n)!}{k!(k+m)!(n-k)!}
\nonumber \\
&& \quad \times[2\psi(k+n+1)-\psi(k+m+1)]
\left(\frac{z\mp1}{2}\right)^{k}
\qquad (0\leqslant m\leqslant n),
\nonumber \\
&&
\label{3.4}
\end{eqnarray}
\begin{eqnarray}
\frac{\mathrm{d}^{m}}{\mathrm{d}z^{m}}[P_{n}(z)\ln(z\pm1)]
&=& \frac{(2m)!}{2^{m}m!}C_{n-m}^{(m+1/2)}(z)\ln(z\pm1)
\nonumber \\
&& -\,(\mp)^{n}(-)^{m}(z\pm1)^{-m}\sum_{k=0}^{m-1}(\pm)^{k}
\frac{(k+n)!(m-k-1)!}{k!(n-k)!}\left(\frac{z\pm1}{2}\right)^{k}
\nonumber \\
&& +\,\frac{(\mp)^{n+m}}{2^{m}}
\sum_{k=0}^{n-m}(\mp)^{k}\frac{(k+n+m)!}{k!(k+m)!(n-m-k)!}
\nonumber \\
&& \quad \times[\psi(k+m+1)-\psi(k+1)]
\left(\frac{z\pm1}{2}\right)^{k}
\qquad (0\leqslant m\leqslant n),
\label{3.5}
\end{eqnarray}
\begin{eqnarray}
\frac{\mathrm{d}^{m}}{\mathrm{d}z^{m}}[P_{n}(z)\ln(z\pm1)]
&=& \frac{(2m)!}{2^{m}m!}C_{n-m}^{(m+1/2)}(z)\ln(z\pm1)
\nonumber \\
&& +\,\frac{(2m)!}{2^{m}m!}[\psi(n+1)-\psi(n-m+1)]C_{n-m}^{(m+1/2)}(z)
\nonumber \\
&& +\,(\mp)^{n}\frac{(n+m)!}{(n-m)!}(z^{2}-1)^{-m/2}
\sum_{k=0}^{m-1}(-)^{k}\frac{2k+1}{(n-k)(k+n+1)}P_{k}^{-m}(\pm z)
\nonumber \\
&& -\,(\mp)^{n+m}\frac{(2m)!}{2^{m}m!}
\sum_{k=0}^{n-m-1}(\mp)^{k}\frac{2k+2m+1}{(n-m-k)(k+n+m+1)}
\nonumber \\
&& \quad \times\left[1-\frac{k!(n+m)!}{(k+2m)!(n-m)!}\right]
C_{k}^{(m+1/2)}(z)
\qquad (0\leqslant m\leqslant n)
\label{3.6}
\end{eqnarray}
and
\begin{equation}
\frac{\mathrm{d}^{m}}{\mathrm{d}z^{m}}[P_{n}(z)\ln(z\pm1)]
=(\mp)^{n}(-)^{m+1}(n+m)!(m-n-1)!(z^{2}-1)^{-m/2}P_{n}^{-m}(\pm z)
\qquad (m>n).
\label{3.7}
\end{equation}
The sextet of formulas (\ref{3.2})--(\ref{3.7}) constitutes the
result of this paper.
\section{Concluding remarks}
\label{IV}
\setcounter{equation}{0}
In some applications, it might be necessary to have explicit
representations of \mbox{$\mathrm{d}^{m}[P_{n}(x)\ln(1\pm
x)]/\mathrm{d}x^{m}$} with $-1\leqslant x\leqslant1$. Such explicit
formulas may be derived from Eqs.\ (\ref{3.2})--(\ref{3.7}) with the
aid of the relationships
\begin{equation}
x+1\pm\mathrm{i}0=1+x,
\qquad 
x-1\pm\mathrm{i}0=\mathrm{e}^{\pm\mathrm{i}\pi}(1-x)
\qquad (-1\leqslant x\leqslant1)
\label{4.1}
\end{equation}
and
\begin{eqnarray}
P_{n}^{\pm m}(x) &=& \mathrm{e}^{\pm\mathrm{i}\pi m/2}
P_{n}^{\pm m}(x+\mathrm{i}0)
=\mathrm{e}^{\mp\mathrm{i}\pi m/2}P_{n}^{\pm m}(x-\mathrm{i}0)
\nonumber \\
&=& \frac{1}{2}\left[\mathrm{e}^{\pm\mathrm{i}\pi m/2}
P_{n}^{\pm m}(x+\mathrm{i}0)+\mathrm{e}^{\mp\mathrm{i}\pi m/2}
P_{n}^{\pm m}(x-\mathrm{i}0)\right]
\qquad (-1\leqslant x\leqslant1).
\label{4.2}
\end{eqnarray}
Since the procedure is straightforward and does not offer any
difficulty, we do not list the resulting expressions here.
%
%

%
\end{document}